\documentclass[11pt]{amsart}

\usepackage{enumerate,url,amssymb, mathrsfs}

\newtheorem{theorem}{Theorem}[section]

\newtheorem{corollary}[theorem]{Corollary}

\theoremstyle{definition}

\newtheorem{question}[theorem]{Question}

\theoremstyle{remark}

\numberwithin{equation}{section}

\newcommand{\abs}[1]{\lvert#1\rvert}

\newcommand{\C}{\mathbb{C}}
\newcommand{\DD}{\mathbb{D}}
\newcommand{\K}{\mathcal{K}}
\newcommand{\R}{\mathbb{R}}

\newcommand{\D}{\partial}

\DeclareMathOperator{\im}{Im}

\def\XXint#1#2#3{{\setbox0=\hbox{$#1{#2#3}{\int}$}
\vcenter{\hbox{$#2#3$}}\kern-.5\wd0}}

\def\le{\leqslant}
\def\ge{\geqslant}

\begin{document}

\title{Area contraction for harmonic automorphisms of the disk}

\author{Ngin-Tee Koh}
\address{Department of Mathematics, Syracuse University, Syracuse,
NY 13244, USA}
\email{nkoh@syr.edu}

\author{Leonid V. Kovalev}
\address{Department of Mathematics, Syracuse University, Syracuse,
NY 13244, USA}
\email{lvkovale@syr.edu}
\thanks{Kovalev was supported by the NSF grant DMS-0913474.}


\date{November 2, 2009}


\begin{abstract}
A harmonic self-homeomorphism of a disk does not increase the area of any concentric disk.
\end{abstract}

\maketitle

\section{Introduction}

The unit disk $\DD=\{z\in\C \colon \abs{z}<1\}$ can be endowed with the hyperbolic metric
\[d\sigma=\frac{\abs{dz}}{1-\abs{z}^2}.\]
The Schwarz-Pick lemma (e.g., ~\cite{Abook}) implies that
any holomorphic map $f\colon \DD\to\DD$ does not increase distances in the hyperbolic metric. This is no longer true
for harmonic maps, which verify the Laplace equation $\D \bar \D f=0$ but not necessarily the Cauchy-Riemann equation $\bar\D f=0$.
The harmonic version of the Schwarz lemma (\cite{He}, see also~\cite{CGH}) states
that any harmonic map $f\colon\DD\to\DD$ with normalization $f(0)=0$ satisfies
\[\abs{f(z)}\le \frac{4}{\pi}\arctan\abs{z},\quad z\in\DD.\]
This inequality is sharp~\cite[p. 77]{Dub}. More precisely, for any $r\in (0,1)$ and any small $\epsilon>0$ there is a
bijective harmonic map $f\colon \DD\to\DD$ such that $f(0)=0$ and
\[f(r)=-f(-r)=\frac{4}{\pi}\arctan r -\epsilon.\]
This map is not a contraction in either Euclidean or hyperbolic metric. With respect to either
metric, the diameter of the disk $\DD_r =\{z\in\C \colon \abs{z}<r\}$ is strictly less than the diameter of $f(\DD_r)$.

In this note we prove that a bijective harmonic map $f\colon \DD\to\DD$ does not increase the area of $\DD_r$ for any $0<r<1$.
We write $\abs{E}$ for the area (i.e., planar Lebesgue measure) of a  set $E$.

\begin{theorem}\label{main}
Let $f\colon \DD\to\DD$ be a bijective harmonic map. Then
\begin{equation}\label{main1}
\abs{f(\DD_r)} \le \abs{\DD_r},\qquad 0< r <1.
\end{equation}
If~\eqref{main1} turns into an equality for some $r\in (0,1)$, then $f$ is an isometry.
\end{theorem}

It should be noted that the class of harmonic automorphisms of $\DD$ is much wider than
the class of holomorphic automorphisms, which consists of M\"obius maps only. Harmonic homeomorphisms of $\DD$
form an interesting and much-studied class of planar maps, see~\cite{CSS,Ka,Pa} or the monograph~\cite{Dub}.
Theorem~\ref{main} is different from most known estimates for harmonic maps in that it remains sharp
when specialized to the holomorphic case.

An immediate consequence of~\eqref{main1} is
$\abs{f(\DD\setminus \DD_r)} \ge \abs{\DD\setminus\DD_r}$. If $f$ is sufficiently smooth,
we can divide by $1-r$ and let $r\to 1$ to obtain the following.

\begin{corollary}\label{cor}
Let $f\colon \DD\to\DD$ be a bijective harmonic map that is continuously differentiable in the closed
disk $\overline{\DD}$. Then
\begin{equation*}
\int_{\abs{z}=1} \abs{\det Df}\,\abs{dz} \ge 2\pi,
\end{equation*}
where $\det Df=\abs{\D f}^2-\abs{\bar\D f}^2$ is the Jacobian determinant of $f$.
\end{corollary}

Corollary~\ref{cor} was proved in a different way in~\cite{IKO2} where it serves as an important part of the
proof of Nitsche's conjecture on the existence of harmonic homeomorphisms between doubly-connected domains.
In fact, Corollary~\ref{cor} is what led us to think that~\eqref{main1} might be true.

If $f\colon \DD\to\DD$ is holomorphic, then~\eqref{main1} holds without the assumption of $f$ being bijective.
Indeed, in this case $f(\DD_r)$ is contained in a hyperbolic disk $D$ of the same hyperbolic radius as $\DD_r$. Since the
density of the hyperbolic metric increases toward the boundary, it follows that the Euclidean radius of $D$ is at most $r$,
which implies~\eqref{main1}.

\begin{question} Does the area comparison~\eqref{main1} hold for general harmonic maps $f\colon \DD\to\DD$?
Does it hold in higher dimensions?
\end{question}

We conclude the introduction by comparing the behavior of $\abs{f(\DD_r)}$ for holomorphic and harmonic maps.
If $f\colon \DD\to \C$ is holomorphic and injective, one can use the power series $f(z)=\sum c_n z^n$ to compute
\[\abs{f(\DD_r)} = \sum_{n=1}^{\infty} n \abs{c_n}^2 r^{2n}.\]
Since the right-hand side is a convex function of $r^2$, it follows that
\begin{equation}\label{general}
\abs{f(\DD_r)}\le r^2 \abs{f(\DD)},
\end{equation}
which includes~\eqref{main1} as a special case. However, ~\eqref{general} fails for harmonic maps.
Indeed, let $f(z)=z+c\bar z^2$ where $0<\abs{c}<1/2$. It is easy to see that $f\colon \DD\to \C$ is harmonic and one-to-one, but
\[\abs{f(\DD_r)} = r^2-2\abs{c}^2 r^4\]
is a strictly concave function of $r^2$. Therefore, $\abs{f(\DD_r)}>r^2 \abs{f(\DD)}$ for $0<r<1$.
This example does not contradict Theorem~\ref{main} since $f(\DD)$ is not a disk.

\section{Preliminaries}

Let $f$ be as in Theorem~\ref{main}. We may assume that $f$ is orientation-preserving; otherwise consider $f(\bar z)$ instead.
In this section we derive an identity that relates the area of $f(\DD_r)$ with the boundary values of $f$, which
exist a.e. in the sense of nontangential limits.

The Poisson kernel for $\DD$ will be denoted $P_{r}(t)$,
\[P_r (t) = \frac{1-r^2}{1-2r\cos t+r^2},\quad 0\le r<1, \ t\in\R.\]
We represent $f$ by the Poisson integral
\begin{equation}\label{Pintegral}
f(r e^{i\theta}) = \frac{\omega}{2\pi} \int_{0}^{2\pi} e^{i\xi(t)}P_r(\theta-t) \, dt,
\end{equation}
where $\xi \colon [0,2\pi)\to [0,2\pi)$ is a nondecreasing function and $\omega$ is a unimodular constant.
By Green's formula we have
\begin{equation*}
\abs{f(\DD_r)}= \frac{1}{2}\int_{0}^{2\pi} \im\big(\overline{f(re^{i\theta})} f_{\theta}(re^{i\theta})\big)\,d\theta,
\end{equation*}
where $f_\theta$ indicates the derivative with respect to $\theta$. Since
\begin{equation*}
f_\theta(r e^{i\theta}) = \frac{\omega}{2\pi} \int_{0}^{2\pi} e^{i\xi(t)}P_r'(\theta-t)\, dt,
\end{equation*}
it follows that
\begin{equation}\label{ar2}
\overline{f(re^{i\theta})} f_{\theta}(re^{i\theta}) =
\frac{1}{4\pi^2}\int_{0}^{2\pi} \int_{0}^{2\pi} e^{-i\xi(t)} e^{i\xi(s)} P_r(\theta-t) P_r'(\theta-s) \, dt\,ds.
\end{equation}
Integrating~\eqref{ar2} with respect to $\theta$ and reversing the order of integration, we find
\begin{equation}\label{ar3}
\abs{f(\DD_r)}=
\frac{1}{4\pi}\int_{0}^{2\pi} \int_{0}^{2\pi} \K_r \, \sin(\xi(s)-\xi(t)) \, dt\,ds
\end{equation}
where $\K_r$ is a function of $r$, $s$, and $t$,
\[
\K_r = \frac{1}{2\pi} \int_{0}^{2\pi} P_r(\theta-t) P_r'(\theta-s)
\]

Recall that the Poisson kernel has the semigroup property~\cite[p.62]{St},
\begin{equation}\label{Pconv}
P_{r \sigma}(t) = \frac{1}{2\pi} \int_{0}^{2\pi} P_r(s)P_\sigma(t-s) \, ds,\quad 0\le r,\sigma<1.
\end{equation}
We will only use~\eqref{Pconv} with $\sigma=r$. Differentiation with respect to $t$ yields
\begin{equation}\label{Pconv2}
\frac{1}{2\pi} \int_{0}^{2\pi} P_r(s)P'_r(t-s)\, ds = P_{r^2}'(t)
= -\frac{2r^2 (1-r^4) \sin t}{(1-2r^2\cos t+r^4)^2}.
\end{equation}
Identity~\eqref{Pconv2} provides an explicit formula for $\K_r$,
\begin{equation}\label{kernel1}
\K_r = \K_r(s-t)= \frac{2\rho^2 (1-\rho^4) \sin (s-t)}{(1-2\rho^2\cos (s-t)+\rho^4)^2}.
\end{equation}
Now we can rewrite~\eqref{ar3} as
\begin{equation}\label{ar4}
\abs{f(\DD_r)}=
\frac{1}{4\pi}\int_{0}^{2\pi} \int_{0}^{2\pi} \K_r(s-t)\, \sin(\xi(s)-\xi(t))\, dt\,ds.
\end{equation}
In the next section we will estimate~\eqref{ar4} from above.

\section{Proof of Theorem~\ref{main}}

We continue to use the Poisson representation~\eqref{Pintegral}.
The function $\xi$, originally defined on $[0,2\pi)$, can be extended
to $\R$ so that $\xi(t+2\pi)=\xi(t)+2\pi$ for all $t\in\R$. By~\eqref{ar4} we have
\begin{equation}\label{ar5}
\abs{f(\DD_r)}=
\frac{1}{4\pi}\int_{0}^{2\pi} \int_{0}^{2\pi} \K_r(s-t)\, \sin(\xi(s)-\xi(t))\, dt\,ds.
\end{equation}
When $f$ is the identity map,~\eqref{ar5} tells us that
\[
\frac{1}{4\pi}\int_{0}^{2\pi} \int_{0}^{2\pi} \K_r(s-t)\, \sin(s-t)\, dt\,ds  = \abs{\DD_r}.
\]
The desired inequality $\abs{f(\DD_r)}\le \abs{\DD_r}$ now takes the form
\begin{equation}\label{ar7}
\int_{0}^{2\pi} \int_{0}^{2\pi} \K_r(s-t)\, \big\{\sin(s-t)-\sin(\xi(s)-\xi(t))\big\}\, dt\,ds \ge 0.
\end{equation}
Neither the kernel $\K_r$, which is defined by~\eqref{kernel1}, nor the other factor in the integrand are nonnegative.
We will have to transform the integral in~\eqref{ar7} before effective pointwise estimates can be made.
It will be convenient to use the notation
\begin{equation}\label{alga}
\alpha=s-t, \quad \text{and  }\  \gamma=\gamma(\alpha,t)=\xi(\alpha+t)-\xi(t),
\end{equation}
so that the integral in~\eqref{ar7} becomes
\[
\int_{0}^{2\pi} \int_{-\pi-t}^{\pi-t} \K_r(\alpha)\, (\sin \alpha-\sin \gamma )\, d\alpha\,dt.
\]
Since the integrand is $2\pi$-periodic with respect to $\alpha$, our goal can be equivalently stated as
\begin{equation}\label{ar17}
\int_{0}^{2\pi} \int_{0}^{2\pi} \K_r(\alpha)\, (\sin \alpha-\sin \gamma)\, d\alpha\,dt \ge 0.
\end{equation}
Note that $\gamma\in [0,2\pi]$ for all $\alpha,t\in [0,2\pi]$.

\textbf{Step 1.} We claim that
\begin{equation}\label{s1t0}
\int_{0}^{2\pi} \int_{0}^{2\pi} \K_r(\alpha)(\gamma  - \alpha)\cos\alpha \,d\alpha\,dt =0.
\end{equation}
Indeed, the function $\zeta(t):=\xi(t)-t$ is $2\pi$-periodic, which implies
\begin{equation}\label{s1t1}
\int_{0}^{2\pi} \{\zeta(\alpha+t)-\zeta(t)\}\,dt=0
\end{equation}
for every $\alpha\in\R$. Multiplying~\eqref{s1t1} by $\K_r(\alpha)\cos \alpha$ and integrating over $\alpha\in [0,2\pi]$, we obtain
\[
\int_{0}^{2\pi} \int_{0}^{2\pi}\K_r(\alpha)\{\zeta(\alpha+t)-\zeta(t)\}\cos \alpha \,d\alpha\,dt=0
\]
It remains to note that $\zeta(\alpha+t)-\zeta(t)=\gamma - \alpha$, completing the proof of~\eqref{s1t0}.

We take advantage of~\eqref{s1t0} by adding it to~\eqref{ar17}, which reduces our task to proving that
\begin{equation}\label{ar8}
\int_{0}^{2\pi} \int_{0}^{2\pi} \K_r(\alpha)\left\{\sin\alpha + (\gamma - \alpha)\cos\alpha - \sin\gamma\right\}\,d\alpha\,dt \ge 0.
\end{equation}

\textbf{Step 2.} Let us now consider the function
\begin{equation}\label{tangent}
H(\alpha,\beta):=\sin\alpha + (\beta - \alpha)\cos\alpha - \sin\beta, \quad (\alpha,\beta)\in [0,2\pi]\times [0,2\pi]
\end{equation}
which appears in~\eqref{ar8}. It has a simple geometric interpretation in terms of the graph of the sine function $y = \sin x$.
Indeed, the tangent line to this graph at $x=\alpha$ has equation $y=\sin\alpha + (x-\alpha)\cos\alpha$.
The quantity $H(\alpha,\beta)$ represents the difference in the $y$-values of the tangent line and the graph at $x=\beta$.
Since the sine curve is strictly concave on $[0,\pi]$, it follows that
\begin{equation}\label{hpos1}
H(\alpha,\beta)\ge 0,\qquad 0\le \alpha,\beta\le \pi,
\end{equation}
with equality only when $\alpha=\beta$. The upper bound on $\beta$ in~\eqref{hpos1} can be weakened
to $\beta\le 2\pi-\alpha$ thanks to the monotonicity with respect to $\beta$,
\[
\frac{\D H}{\D \beta}=\cos\alpha - \cos\beta \ge 0,\quad 0\le\alpha\le \pi, \  \alpha\le \beta\le 2\pi-\alpha.
\]
Note that the product $\K_r(\alpha)H(\alpha,\beta)$ is invariant under the central symmetry of the square
$[0,2\pi]\times [0,2\pi]$, i.e., the transformation $(\alpha,\beta)\mapsto (2\pi-\alpha,2\pi-\beta)$. Hence
\begin{equation}\label{pointwise}
\K_r(\alpha)H(\alpha,\beta)\ge 0, \qquad (\alpha,\beta) \in \big([0,2\pi]\times [0,2\pi]\big) \setminus (T_1\cup T_2)
\end{equation}
where
\begin{align*}
T_1&=\{(\alpha,\beta)\colon 0< \alpha < \pi, \ 2\pi-\alpha < \beta \le 2\pi\}; \\
T_2&=\{(\alpha,\beta)\colon \pi < \alpha < 2\pi, \ 0 \le \beta < 2\pi-\alpha\}.
\end{align*}
Within the triangles $T_1$ and $T_2$ the product $\K_r(\alpha)H(\alpha,\beta)$ may be negative.
However, for all $(\alpha,\beta)\in [0,2\pi]\times [0,2\pi]$ the following holds.
\begin{equation}\label{symsum}
\K_r(\alpha)H(\alpha,\beta) + \K_r(2\pi-\alpha)H(2\pi-\alpha,\beta) = 2\K_r(\alpha)H(\alpha,\pi) \ge 0,
\end{equation}
where the last inequality follows from~\eqref{pointwise}. We will use~\eqref{symsum} to control
the contribution of triangles $T_1$ and $T_2$ to the integral~\eqref{ar8}.

\textbf{Step 3.} For each fixed $t$ the function $\alpha\mapsto \gamma(\alpha,t)$ defined by~\eqref{alga}
is nondecreasing and it maps the interval $[0,2\pi]$ onto itself. Thus,
inequality~\eqref{ar8} will follow once we show that for any nondecreasing function
$\Gamma \colon [0,2\pi]\to [0,2\pi]$
\begin{equation}\label{ar10}
\int_{0}^{2\pi} \K_r(\alpha) H(\alpha,\Gamma(\alpha))\, d\alpha \ge 0.
\end{equation}
The integral in~\eqref{ar10} remains unchanged if we replace $\Gamma(\alpha)$ with
$\widetilde{\Gamma}(\alpha)=2\pi-\Gamma(2\pi-\alpha)$. Thus we lose no generality in assuming that
$\Gamma(\pi)\le \pi$. By virtue of~\eqref{pointwise} the integrand in~\eqref{ar10} is nonnegative outside of the interval $[\pi,\alpha_0]$, where
\[
\alpha_0=\sup\{\alpha\in [\pi,2\pi]\colon \alpha+\Gamma(\alpha)\le 2\pi\}
\]
We claim that
\begin{equation}\label{comp1}
\K_r(\alpha) H(\alpha,\Gamma(\alpha)) \ge \K_r(\alpha) H(\alpha,\Gamma(\pi)), \quad
2\pi-\alpha_0< \alpha < \alpha_0.
\end{equation}
Indeed, the inequality
\[
\frac{\D H}{\D \beta}=\cos\alpha - \cos\beta \le 0,\quad \abs{\alpha-\pi}\le \abs{\beta-\pi}\le \pi,
\]
implies
\begin{equation}\label{hmot}
H(\alpha,\beta_1) \ge H(\alpha,\beta_2),\quad
0\le \beta_1\le \beta_2 \le \min(\alpha,2\pi-\alpha).
\end{equation}
To see that~\eqref{hmot} applies in our situation, note that
$\Gamma(\alpha)\le 2\pi-\alpha_0$ for $\alpha<\alpha_0$.
Inequality~\eqref{hmot} yields
\begin{equation}\label{comp2}
\begin{split}
H(\alpha,\Gamma(\alpha)) &\le H(\alpha,\Gamma(\pi)), \quad \pi\le \alpha < \alpha_0; \\
H(\alpha,\Gamma(\alpha)) &\ge H(\alpha,\Gamma(\pi)), \quad 2\pi-\alpha_0 <  \alpha \le \pi.
\end{split}
\end{equation}
Multiplying~\eqref{comp2} by $\K_r(\alpha)$, we arrive at~\eqref{comp1}.

Finally, we combine~\eqref{pointwise}, \eqref{comp1}, and~\eqref{symsum} to obtain
\begin{equation}\label{punchline} \begin{split}
\int_{0}^{2\pi} \K_r(\alpha) H(\alpha,\Gamma(\alpha))\, d\alpha
&\ge \int_{2\pi-\alpha_0}^{\alpha_0} \K_r(\alpha) H(\alpha,\Gamma(\alpha))\, d\alpha \\
&\ge \int_{2\pi-\alpha_0}^{\alpha_0} \K_r(\alpha) H(\alpha,\Gamma(\pi))\, d\alpha  \\
&= 2 \int_{\pi}^{\alpha_0} \K_r(\alpha) H(\alpha,\pi)\, d\alpha \ge 0,
\end{split}
\end{equation}
completing the proof of~\eqref{ar8}.

\textbf{Step 4.} It remains to prove the equality statement in Theorem~\ref{main}.
Suppose that $\Gamma \colon [0,2\pi]\to [0,2\pi]$ is a nondecreasing function such that $\Gamma(\pi)\le \pi$, and
equality holds everywhere in~\eqref{punchline}. Returning to
the geometric interpretation of $H(\alpha,\gamma)$ in~\eqref{tangent}, we note that
\[\K_r(\alpha) H(\alpha,\pi)>0,\quad 0<\abs{\alpha-\pi}<\pi.\]
This forces $\alpha_0=\pi$, which by definition of $\alpha_0$ implies
\begin{equation}\label{equal1}
\K_r(\alpha) H(\alpha,\Gamma(\alpha))\ge 0,\quad 0\le\alpha\le 2\pi.
\end{equation}
Hence, ~\eqref{equal1} must turn into an equality for almost all $\alpha\in [0,2\pi]$.
In view of~\eqref{hpos1} and of the monotonicity of $\Gamma$ this is only possible if
$\Gamma(\alpha)=\alpha$ for all $\alpha\in [0,2\pi]$.

If $\abs{f(\DD_r)}=\abs{\DD_r}$, then equality holds in~\eqref{ar8}. Then for almost all $t\in [0,2\pi]$
the function $\Gamma(\alpha)=\xi(\alpha+t)-\xi(t)$, or its reflection
$\widetilde{\Gamma}(\alpha)=2\pi-\Gamma(2\pi-\alpha)$,  turns~\eqref{punchline} into an equality.
Hence $\xi(\alpha+t)-\xi(t)=\alpha$ for almost all $t\in [0,2\pi]$ and all $\alpha\in [0,2\pi]$.
Thus $\xi$ is the identity function and $f\colon\DD\to\DD$ is an isometry. Theorem~\ref{main} is proved.

\section*{Acknowledgements}
We thank Tadeusz Iwaniec and Jani Onninen for valuable discussions on the subject of this paper.

\bibliographystyle{amsplain}

\end{document}